\documentclass[11pt]{amsart}
\usepackage[T1]{fontenc}
\usepackage{amsmath,amssymb,amsthm,mathtools}
\usepackage{enumitem}
\usepackage{xcolor}
\usepackage[colorlinks=true,citecolor=black,linkcolor=black,urlcolor=black]{hyperref}
\usepackage[margin=1.12in]{geometry}
\numberwithin{equation}{section}
\newtheorem{theorem}{Theorem}[section]
\newtheorem{lemma}[theorem]{Lemma}
\newtheorem{corollary}[theorem]{Corollary}

\theoremstyle{definition}
\theoremstyle{remark}
\newcommand{\Hcal}{\mathcal H}\newcommand{\Fcal}{\mathcal F}
\newcommand{\Acal}{\mathcal A}\newcommand{\Bcal}{\mathcal B}
\newcommand{\Rcal}{\mathcal R}\newcommand{\miss}{\operatorname{miss}}

\newcommand{\e}{\mathrm e}
\title[Quadratic bounds for intersecting hypergraphs]{\texorpdfstring{Quadratic bounds for intersecting hypergraphs: an affirmative answer to a question of Duke and Erd\H{o}s}{Quadratic bounds for intersecting hypergraphs: an affirmative answer to a question of Duke and Erdős}}
\author{Guoqing Wang}\author{Linlin Wang}\author{Jun Yue}
\date{}
\begin{document}
\begin{abstract}
Let $m(k,b)$ denote the minimum number of edges in an intersecting
$k$-uniform hypergraph in which every $(k-1)$-set of vertices is
disjoint from at least $b$ edges, with $m(k,b)=+\infty$ if no such
hypergraph exists. Duke and Erd\H{o}s
[Congressus Numerantium (1980)]
asked whether $m(k,k)=o(k^4)$ is possible.
In this paper, we answer this question affirmatively and determine
the order of growth of $m(k,k)$ by proving that
$k^2+1\leq m(k,k)\leq Ck^2$ for all sufficiently large $k$,
where $C$ is an absolute constant. More generally, we show that
$m(k,b)=\Theta(kb)$ for all sufficiently large odd $k$ and every
integer $b$ with $k\leq b\leq (k-1)(k-3)/2$,
where the implied constants are independent of $k$ and $b$.
\end{abstract}

\maketitle
\noindent\textbf{Keywords.} intersecting hypergraphs; covering number; diversity; probabilistic method.\\
\textbf{2020 Mathematics Subject Classification.} 05D05, 05C65.

\section{Introduction}
All hypergraphs considered in this paper are finite and have no repeated edges.
We identify a hypergraph $\Hcal$ with its edge family and write $|\Hcal|$
for its number of edges. A hypergraph is $k$-uniform if every edge has
size $k$, and intersecting if every two edges have nonempty intersection.
A cover is a vertex set meeting every edge, and the covering number
$\tau(\Hcal)$ is the minimum size of a cover. For $S\subseteq V(\Hcal)$,
write $\miss_{\Hcal}(S)=|\{E\in\Hcal:E\cap S=\varnothing\}|$.
The $t$-diversity of $\Hcal$ is
\[
\gamma_t(\Hcal)
=\min_{S\in\binom{V(\Hcal)}{t}}\miss_{\Hcal}(S),
\qquad 0\leq t\leq |V(\Hcal)|.
\]
For integers $k\geq2$ and $b\geq1$, define
\[
m(k,b)=\min\{|\Hcal|:\Hcal\text{ is intersecting and $k$-uniform, }
\gamma_{k-1}(\Hcal)\geq b\},
\]
with $m(k,b)=+\infty$ if no such hypergraph exists.
We abbreviate $m(k,k)$ to $m(k)$ and write $g(k)=m(k,1)$.
Since every edge of a nonempty intersecting $k$-uniform hypergraph
is a cover, the condition $\gamma_{k-1}(\Hcal)\geq1$ is equivalent
to $\tau(\Hcal)=k$. Thus $g(k)$ is the minimum number of edges
in an intersecting $k$-uniform hypergraph whose covering number
equals its uniformity.

The study of $m(k)$ goes back to Duke and Erd\H{o}s~\cite{DE80},
who considered Ramsey numbers involving matchings and almost complete
hypergraphs. They investigated the assertion that a $k$-uniform
hypergraph in which every $(k-1)$-set misses at least $k$ edges
must contain two disjoint edges. They proved this assertion for
$k=2,3$, but constructed counterexamples for larger uniformities.
One construction uses three disjoint projective planes, with edges
formed by taking the union of a line from each of two distinct planes.
For certain even values of $k$, it yields an intersecting counterexample
with $O(k^4)$ edges. This led Duke and Erd\H{o}s to ask for the
minimum number of edges in such a counterexample and, in particular,
whether $m(k)=o(k^4)$ is possible~\cite[Question~2, p.~374]{DE80}.
Their results for $k=2,3$ give $m(2)=m(3)=+\infty$.

The case of avoidance multiplicity one is the classical covering-number
problem studied by Erd\H{o}s and Lov\'asz~\cite{EL75}.
Kahn~\cite{Kahn94} proved that $g(k)=O(k)$, determining its order
of growth. More recent work has addressed refinements and additional
constraints. Buci\'c, Jain, and Sivashankar~\cite{BJS25} studied
the minimum number of edges when the number of vertices is prescribed,
while Sivashankar~\cite{Siv26} established $g(k)\geq3k-4$ and the
stronger asymptotic bound
$g(k)\geq((41-\sqrt{19})/12-o(1))k$.
These results require every $(k-1)$-set to miss at least one edge;
the parameter $m(k,b)$ extends this requirement to at least $b$ edges.

The same avoidance condition can also be expressed in terms of
diversity. Frankl and Wang~\cite{FW25} studied upper bounds on
diversity and employed wreath products of intersecting families.
Here we seek the minimum number of edges needed to attain a prescribed
diversity. For even uniformities, the three-block construction of
Duke and Erd\H{o}s~\cite[p.~373]{DE80} can be applied to an
intersecting $r$-uniform family with covering number $r$.
The avoidance estimate in Section~\ref{sec:lower} gives
$m(2r)\leq3g(r)^2$ for $r\geq3$, and Kahn's theorem therefore
implies $m(2r)=O(r^2)$. Hence classical results already yield
a quadratic bound for even uniformities. To obtain the same order
for all sufficiently large uniformities, it remains to treat
the odd case.

In this paper, we establish such a bound, thereby answering
the question of Duke and Erd\H{o}s affirmatively and determining
the order of growth of $m(k)$.
\begin{theorem}\label{thm:main}
There are absolute constants $C>0$ and $k_0$ such that
$k^2+1\leq m(k)\leq Ck^2$ for every $k\geq k_0$.
In particular, $m(k)=\Theta(k^2)$.
\end{theorem}

The lower bound follows from the general inequality
$m(k,b)\geq kb+1$, valid for all $k\geq2$ and $b\geq1$.
We also show that this inequality has the correct order over
a range of larger avoidance multiplicities. More precisely,
for all sufficiently large odd $k$,
\[
m(k,b)=\Theta(kb)
\qquad\text{whenever}\qquad
k\leq b\leq\frac{(k-1)(k-3)}2,
\]
where the implied constants are independent of $k$ and $b$;
see Corollary~\ref{cor:more}.

To prove the odd-uniformity bound, we use a cyclic construction
based on pairs of cross-intersecting families, in which every
member of one family meets every member of the other.
Starting from Kahn's family, we construct an intersecting
$(2r+1)$-uniform hypergraph with $O(r^3)$ edges such that
every $2r$-set misses at least $2r(r-1)$ edges.
Although there may be many $2r$-sets, they determine only
$\exp(O(r))$ distinct sets of missed edges.
This observation allows a standard sampling argument to retain
$O(r^2)$ edges while ensuring that every $2r$-set still misses
at least $2r+1$ edges. Varying the sampling probability gives
the corresponding bound for larger avoidance multiplicities.

The remainder of the paper is organized as follows.
Section~\ref{sec:lower} establishes the elementary estimates,
states Kahn's theorem, and derives the even-uniformity bound
from the three-block construction.
Section~\ref{sec:odd} gives the odd-uniformity construction
and its sparsification, completing the proof of
Theorem~\ref{thm:main}. The same section then establishes
the extension to larger avoidance multiplicities and
concludes with further questions.

\section{Preliminaries and even uniformities}\label{sec:lower}
We begin with two elementary estimates. The first bounds the
number of edges missed by a vertex set in terms of its size
and the covering number. The second gives the lower bounds
for the extremal functions considered in this paper.

\begin{lemma}\label{lem:miss-cover}
Let $\Fcal$ be an intersecting hypergraph with $\tau(\Fcal)=r$.
For every vertex set $T$ with $|T|=t<r$, we have
$\miss_{\Fcal}(T)\geq2(r-t)-1$.
Moreover, if $\Fcal$ is $r$-uniform and $r\geq3$, then
$|\Fcal|\geq2r$.
\end{lemma}

\begin{proof}
Let $\Fcal_T=\{F\in\Fcal:F\cap T=\varnothing\}$ and
$h=|\Fcal_T|=\miss_{\Fcal}(T)$.
Every cover $C$ of $\Fcal_T$ gives a cover $C\cup T$ of $\Fcal$,
so $r\leq |C\cup T|\leq |C|+t$.
It follows that $\tau(\Fcal_T)\geq r-t$.
Since $\Fcal_T$ is intersecting, it also has a cover of size
at most $\lceil h/2\rceil$: pair its edges, choose a common
vertex from each pair, and, if an edge remains, choose one
of its vertices.
Thus $r-t\leq\lceil h/2\rceil$, which implies
$h\geq2(r-t)-1$.

Now suppose that $\Fcal$ is $r$-uniform and $r\geq3$.
Taking $T=\varnothing$ gives $|\Fcal|\geq2r-1$.
Suppose for a contradiction that equality holds.
If a vertex $v$ belongs to at least three edges, choose three
such edges and partition the remaining $2r-4$ edges into
$r-2$ pairs. A common vertex from each pair, together with $v$,
then covers $\Fcal$ with at most $r-1$ vertices,
contrary to $\tau(\Fcal)=r$.
Hence every vertex belongs to at most two edges.
For a fixed edge $E\in\Fcal$, the sets $E\cap F$,
$F\in\Fcal\setminus\{E\}$, are therefore nonempty and pairwise
disjoint. There are $2r-2$ such sets, so
$|E|\geq2r-2>r$, contradicting the $r$-uniformity of $\Fcal$.
Consequently, $|\Fcal|\geq2r$.
\end{proof}

We next apply the avoidance requirement to the
$(k-1)$-subsets of a fixed edge.

\begin{lemma}\label{lem:lower}
Let $k\geq2$ and $b\geq1$ be integers.
If $\Hcal$ is an intersecting $k$-uniform hypergraph with
$\gamma_{k-1}(\Hcal)\geq b$, then $|\Hcal|\geq kb+1$.
Consequently, $m(k,b)\geq kb+1$ and, in particular,
$m(k)\geq k^2+1$.
\end{lemma}

\begin{proof}
Fix an edge $E\in\Hcal$. For each $v\in E$, let
$\Hcal_v=\{F\in\Hcal:F\cap E=\{v\}\}$.
Since $\Hcal$ is intersecting, $\Hcal_v$ consists precisely
of the edges disjoint from $E\setminus\{v\}$.
The hypothesis therefore gives $|\Hcal_v|\geq b$.
These families are pairwise disjoint, and none contains $E$,
since $|E|=k\geq2$. Hence
\[
|\Hcal|\geq1+\sum_{v\in E}|\Hcal_v|\geq kb+1.
\]
Taking the minimum over all admissible hypergraphs yields
$m(k,b)\geq kb+1$, with the same inequality holding when
$m(k,b)=+\infty$. The bound for $m(k)$ follows by setting $b=k$.
\end{proof}

To obtain the corresponding upper bound for even uniformities,
we combine the three-block construction of Duke and
Erd\H{o}s~\cite[p.~373]{DE80} with the following theorem,
which provides an intersecting base family with linearly
many edges and covering number equal to its uniformity.

\begin{theorem}[Kahn~\cite{Kahn94}]\label{thm:kahn}
There exist an absolute constant $C_0>0$ and an integer $r_0$
such that, for every integer $r\geq r_0$, there is an
intersecting $r$-uniform hypergraph $\Fcal$ satisfying
$\tau(\Fcal)=r$ and $|\Fcal|\leq C_0r$.
\end{theorem}

{
We next describe the three-block construction for an arbitrary
intersecting base family of maximum covering number.

\medskip
\noindent\textbf{The three-block construction.}
Let $r\geq3$, and let $\Fcal$ be an intersecting $r$-uniform
hypergraph with $\tau(\Fcal)=r$. Write $m=|\Fcal|$.
Take three copies $\Fcal_0,\Fcal_1,\Fcal_2$ of $\Fcal$
on pairwise disjoint vertex sets $X_0,X_1,X_2$.
On $X_0\cup X_1\cup X_2$, define
\[
\Hcal(\Fcal)=\{E_i\cup E_j:
E_i\in\Fcal_i,\ E_j\in\Fcal_j,\ 0\leq i<j\leq2\}.
\]
Each edge of $\Hcal(\Fcal)$ meets exactly two blocks, and its
intersections with those blocks uniquely determine the two
base edges. There are three choices of a pair of blocks and
$m^2$ choices of the two base edges for each pair.
Thus $\Hcal(\Fcal)$ is $2r$-uniform and has exactly $3m^2$
distinct edges. It is also intersecting: the pairs of blocks
used by any two edges share a block, and their intersections
with that block are intersecting base edges.

The following lemma establishes the avoidance property
needed to apply this construction to the function $m(2r)$.
\par}

\begin{lemma}\label{lem:even-avoid}
{
Let $r\geq3$, and let $\Fcal$ be an intersecting $r$-uniform
hypergraph with $\tau(\Fcal)=r$.
Then the three-block construction $\Hcal(\Fcal)$ satisfies
$\gamma_{2r-1}\bigl(\Hcal(\Fcal)\bigr)\geq2r$.
}
\end{lemma}

\begin{proof}
{Write $\Hcal=\Hcal(\Fcal)$ and $m=|\Fcal|$.}
Fix $S\subseteq V(\Hcal)$ with $|S|=2r-1$.
For $i\in\{0,1,2\}$, put $s_i=|S\cap X_i|$ and
$b_i=\miss_{\Fcal_i}(S\cap X_i)$.
Then $s_0+s_1+s_2=2r-1$.
An edge $E_i\cup E_j$ avoids $S$ precisely when both
$E_i$ and $E_j$ avoid $S$. Since the three types of edges
are distinct, we have
\begin{equation}\label{eq:old-3-1}
\miss_{\Hcal}(S)=b_0b_1+b_1b_2+b_2b_0.
\end{equation}

Suppose first that $s_i<r$ for every $i$.
By Lemma~\ref{lem:miss-cover},
$b_i\geq2(r-s_i)-1\geq1$, and hence
\[
B:=b_0+b_1+b_2
\geq6r-2(s_0+s_1+s_2)-3=2r-1.
\]
Writing $b_i=1+x_i$ with $x_i\geq0$ gives
\[
b_0b_1+b_1b_2+b_2b_0
=2B-3+x_0x_1+x_1x_2+x_2x_0
\geq2B-3.
\]
Therefore $\miss_{\Hcal}(S)\geq4r-5\geq2r$, since $r\geq3$.

It remains to consider the case $s_i\geq r$ for some $i$.
Let $\{j,\ell\}=\{0,1,2\}\setminus\{i\}$.
Then $s_j+s_\ell\leq r-1$, so both $s_j$ and $s_\ell$
are smaller than $r$.
If one of them is zero, say $s_j=0$, then $b_j=m\geq2r$
by Lemma~\ref{lem:miss-cover}, while $b_\ell\geq1$.
It follows from \eqref{eq:old-3-1} that
$\miss_{\Hcal}(S)\geq b_jb_\ell\geq2r$.

We may therefore assume that $s_j,s_\ell\geq1$.
Since their sum is at most $r-1$, each is at most $r-2$.
Set $x=r-s_j$ and $y=r-s_\ell$.
Then $x,y\geq2$ and $x+y\geq r+1$.
Applying Lemma~\ref{lem:miss-cover} in the two blocks yields
\[
\begin{aligned}
b_jb_\ell
&\geq(2x-1)(2y-1)\\
&=4(x-2)(y-2)+6(x+y)-15\\
&\geq6r-9
\geq2r,
\end{aligned}
\]
where the last inequality uses $r\geq3$.
Again, \eqref{eq:old-3-1} gives
$\miss_{\Hcal}(S)\geq b_jb_\ell\geq2r$.
Since $S$ was arbitrary, the lemma follows.
\end{proof}

{
We now apply the construction to a base family supplied by
Theorem~\ref{thm:kahn}. Its linear bound on the number of base
edges, together with Lemma~\ref{lem:even-avoid}, yields the
quadratic bound for even uniformities.
\par}

\begin{theorem}\label{thm:even}
There is an absolute constant $C_1>0$ such that
$m(2r)\leq C_1r^2$ for all sufficiently large $r$.
In particular, one may take $C_1=3C_0^2$, where $C_0$
is the constant in Theorem~\ref{thm:kahn}.
\end{theorem}

\begin{proof}
For all sufficiently large $r$, Theorem~\ref{thm:kahn}
provides an intersecting $r$-uniform hypergraph $\Fcal$
with $\tau(\Fcal)=r$ and $m:=|\Fcal|\leq C_0r$.
{Let $\Hcal=\Hcal(\Fcal)$ be the hypergraph obtained
from the three-block construction.}
The resulting hypergraph $\Hcal$ is intersecting and
$2r$-uniform, with exactly $3m^2$ edges.
By Lemma~\ref{lem:even-avoid}, it also satisfies
$\gamma_{2r-1}(\Hcal)\geq2r$.
Thus $\Hcal$ is admissible in the definition of $m(2r)$,
and hence
\[
m(2r)\leq|\Hcal|=3m^2\leq3C_0^2r^2.
\]
The conclusion follows with $C_1=3C_0^2$.
\end{proof}

\section{Odd uniformities and sparsification}\label{sec:odd}
We now turn to odd uniformities. The construction below uses
two cross-intersecting families of consecutive uniformities
to produce an intersecting $(2r+1)$-uniform hypergraph.
We first establish its avoidance properties and then apply
a sampling argument to reduce its number of edges.

\medskip
\noindent{\textbf{The cyclic construction.}}
Let $r\geq3$, and let $\Fcal$ be an intersecting $r$-uniform
hypergraph with $\tau(\Fcal)=r$ and $m=|\Fcal|$.
Take a set $Z$ of $r$ new vertices, disjoint from $V(\Fcal)$,
and define
\[
\Acal=\Fcal\cup\{Z\},\qquad
\Bcal=\{F\cup\{z\}:F\in\Fcal,\ z\in Z\}.
\]
Then $\Acal$ is $r$-uniform with $m+1$ edges, and $\Bcal$
is $(r+1)$-uniform with $mr$ edges.
These families are cross-intersecting: every member of
$\Fcal$ meets the $\Fcal$-part of every member of $\Bcal$,
while $Z$ meets every member of $\Bcal$ in its added vertex.
Moreover, $\Bcal$ is intersecting because $\Fcal$ is intersecting.

Let $X=V(\Fcal)\cup Z$, and take three copies
$(\Acal_i,\Bcal_i)$ of $(\Acal,\Bcal)$ on pairwise disjoint
vertex blocks $X_i$, $i=0,1,2$.
{Define the cyclic construction $\Hcal_{\mathrm{odd}}(\Fcal)$ by}
\begin{equation}\label{eq:old-4-1}
{\Hcal_{\mathrm{odd}}(\Fcal)}=(\Acal_0\times\Bcal_1)\cup
      (\Acal_1\times\Bcal_2)\cup
      (\Acal_2\times\Bcal_0),
\end{equation}
where $\Acal_i\times\Bcal_j
=\{A\cup B:A\in\Acal_i,\ B\in\Bcal_j\}$.
{For a fixed base family $\Fcal$, write
$\Hcal=\Hcal_{\mathrm{odd}}(\Fcal)$.
The auxiliary set $Z$ and the three disjoint copies are fixed
as part of the construction; their particular labels do not
affect any of the properties below.}
Each edge has size $2r+1$ and meets exactly two blocks.
The pair of blocks determines its type, and its intersections
with those blocks determine its two constituent edges.
Consequently, $\Hcal$ has exactly
\begin{equation}\label{eq:old-4-2}
M=3m(m+1)r
\end{equation}
distinct edges.
Two edges of the same type meet in their $\Bcal$ block.
For two edges of different types, their intersections with
the common block belong to $\Acal_i$ and $\Bcal_i$,
respectively, and hence intersect.
Thus $\Hcal$ is intersecting.

To estimate the number of edges of $\Hcal$ missed by a
vertex set, we first record the corresponding counts
within a single block.
For $T\subseteq X$, put $s=|T|$, $z=|T\cap Z|$, and
$h=\miss_{\Fcal}(T\cap V(\Fcal))$.
The edge $Z$ avoids $T$ precisely when $z=0$,
and each of the $h$ base edges avoiding $T$ can be extended
to a member of $\Bcal$ avoiding $T$ in exactly $r-z$ ways.
Therefore, writing $a(T)=\miss_{\Acal}(T)$ and
$b(T)=\miss_{\Bcal}(T)$, we obtain
\begin{equation}\label{eq:old-4-3}
a(T)=h+\mathbf{1}_{\{z=0\}},\qquad
b(T)=h(r-z).
\end{equation}

The following estimates give lower bounds on the avoidance
counts in a block containing at most $r$ selected vertices.

\begin{lemma}\label{lem:local-odd}
{Let $r\geq3$, and let $\Fcal$ be an intersecting
$r$-uniform hypergraph with $\tau(\Fcal)=r$.
Let $\Acal$, $\Bcal$, and $X$ be as in the cyclic construction,
and define $a(T)=\miss_{\Acal}(T)$ and
$b(T)=\miss_{\Bcal}(T)$ for $T\subseteq X$.
Put $s=|T|$.}
If $s<r$, then
\begin{equation}\label{eq:old-4-4}
a(T)\geq2(r-s),\qquad
b(T)\geq r\bigl(2(r-s)-1\bigr).
\end{equation}
If $s=r$, then $a(T)\geq1$.
\end{lemma}
\begin{proof}
Fix $T\subseteq X$ with $|T|=s\leq r$, and put
$z=|T\cap Z|$ and
$h=\miss_{\Fcal}(T\cap V(\Fcal))$.
Then $|T\cap V(\Fcal)|=s-z$, and
\eqref{eq:old-4-3} gives
$a(T)=h+\mathbf{1}_{\{z=0\}}$ and $b(T)=h(r-z)$.

Suppose first that $s<r$, and put $d=r-s\geq1$.
Since $s-z<r$, Lemma~\ref{lem:miss-cover} gives
\[
h\geq2\bigl(r-(s-z)\bigr)-1=2(d+z)-1.
\]
If $z=0$, then $a(T)=h+1\geq2d$.
If $z\geq1$, then $a(T)=h\geq2d+1$.
Thus $a(T)\geq2d=2(r-s)$ in either case.

For the bound on $b(T)$, note that $0\leq z\leq s<r$,
so $r-z>0$. Using \eqref{eq:old-4-3} and $d=r-s$,
we obtain
\[
\begin{aligned}
b(T)
&=h(r-z)\\
&\geq\bigl(2(d+z)-1\bigr)(r-z)\\
&=r(2d-1)+z\bigl(2(s-z)+1\bigr)\\
&\geq r(2d-1).
\end{aligned}
\]
The last inequality follows from $z\geq0$ and $s-z\geq0$.
Since $d=r-s$, this proves \eqref{eq:old-4-4}.

Finally, suppose that $s=r$.
If $z=0$, then $Z\in\Acal$ is disjoint from $T$,
so $a(T)\geq1$.
If $z\geq1$, then
$|T\cap V(\Fcal)|=r-z<r=\tau(\Fcal)$.
Thus $T\cap V(\Fcal)$ is not a cover of $\Fcal$,
and there exists $F\in\Fcal$ disjoint from
$T\cap V(\Fcal)$.
Since $F\subseteq V(\Fcal)$ and $V(\Fcal)\cap Z=\varnothing$,
the edge $F$ is also disjoint from $T\cap Z$,
and hence from $T$.
As $F\in\Acal$, we again obtain $a(T)\geq1$.
\end{proof}

We now combine the estimates within individual blocks to obtain
a quadratic lower bound on the number of edges missed by
every $2r$-set.

\begin{lemma}\label{lem:odd-avoid}
{Let $r\geq3$, and let $\Fcal$ be an intersecting
$r$-uniform hypergraph with $\tau(\Fcal)=r$.
Then the cyclic construction $\Hcal_{\mathrm{odd}}(\Fcal)$
defined in \eqref{eq:old-4-1} satisfies}
\begin{equation}\label{eq:old-4-5}
{\gamma_{2r}\bigl(\Hcal_{\mathrm{odd}}(\Fcal)\bigr)}\geq D,
\qquad D:=2r(r-1).
\end{equation}
\end{lemma}

\begin{proof}
{Write $\Hcal=\Hcal_{\mathrm{odd}}(\Fcal)$ and
$m=|\Fcal|$, and use the blocks $X_i$ from its construction.
The functions $a$ and $b$ on each copy $X_i$ are understood
via its fixed identification with $X$.}
Fix $S\subseteq V(\Hcal)$ with $|S|=2r$.
For $i\in\{0,1,2\}$, put $s_i=|S\cap X_i|$,
$a_i=a(S\cap X_i)$, and $b_i=b(S\cap X_i)$.
Then $s_0+s_1+s_2=2r$, and the definition of $\Hcal$ gives
\begin{equation}\label{eq:old-4-6}
\miss_{\Hcal}(S)=a_0b_1+a_1b_2+a_2b_0.
\end{equation}
Throughout the proof, subscripts are taken modulo three.
We distinguish three cases according to the number of
indices $i$ for which $s_i\geq r$.

Suppose first that $s_i<r$ for every $i$.
Set $d_i=r-s_i$. Then $d_i\geq1$ and
$d_0+d_1+d_2=r$.
By Lemma~\ref{lem:local-odd},
\[
\miss_{\Hcal}(S)
\geq\sum_{i=0}^2 2d_i\,r(2d_{i+1}-1)
=2r\left(2\sum_{i=0}^2d_id_{i+1}-r\right).
\]
Write $d_i=1+x_i$, where $x_i\geq0$ and
$x_0+x_1+x_2=r-3$.
Expanding the sum yields
\[
\sum_{i=0}^2d_id_{i+1}
=3+2\sum_{i=0}^2x_i+\sum_{i=0}^2x_ix_{i+1}
\geq2r-3.
\]
Consequently,
$\miss_{\Hcal}(S)\geq6r(r-2)\geq2r(r-1)=D$,
where the last inequality holds for $r\geq3$.

Next suppose that $s_i\geq r$ for exactly one index $i$.
Then $s_{i+1},s_{i+2}<r$ and
$s_{i+1}+s_{i+2}\leq r$.
Put $x=r-s_{i+1}$ and $y=r-s_{i+2}$.
Thus $x,y\geq1$ and $x+y\geq r$.
The term $a_{i+1}b_{i+2}$ occurs in
\eqref{eq:old-4-6}, and Lemma~\ref{lem:local-odd} gives
\[
a_{i+1}b_{i+2}
\geq2x\,r(2y-1)
\geq2rxy
\geq2r(x+y-1)
\geq2r(r-1)=D.
\]
Here we used $2y-1\geq y$ and
$xy\geq x+y-1$, the latter following from
$(x-1)(y-1)\geq0$.
Hence $\miss_{\Hcal}(S)\geq D$ in this case as well.

Finally, suppose that at least two of the $s_i$ are at least $r$.
Since their sum is $2r$, two block sizes equal $r$ and
the third equals zero.
Let $X_j$ be the empty block.
Then $b_j=|\Bcal_j|=mr$, while $s_{j-1}=r$ gives
$a_{j-1}\geq1$ by Lemma~\ref{lem:local-odd}.
Using the term $a_{j-1}b_j$ in \eqref{eq:old-4-6}
and the bound $m\geq2r$ from Lemma~\ref{lem:miss-cover},
we obtain
$\miss_{\Hcal}(S)\geq mr\geq2r^2>D$.
Thus every $2r$-set misses at least $D$ edges,
proving the lemma.
\end{proof}

Choosing $\Fcal$ as in Theorem~\ref{thm:kahn} gives
$m=O(r)$, so the construction has $O(r^3)$ edges.
We next use its quadratic avoidance bound to select
a subfamily with only $O(r^2)$ edges in which every
$2r$-set still misses at least $2r+1$ edges.

To carry out this sparsification, we regard the edge set
of the construction as a ground set. Each $2r$-set of
vertices determines a subset of this ground set consisting
of the edges it misses, and our selection must retain at
least $2r+1$ elements of every such subset.
The following elementary sampling lemma provides a
selection whose size is controlled by the number of
distinct subsets and their minimum size.

\begin{lemma}\label{lem:selection}
Let $\Omega$ be a finite set of size $M$, let $b\geq2$
be an integer, and let $L\geq0$ and $D\geq b$ be real numbers.
Suppose that $\Rcal$ is a family of at most $\e^L$ subsets
of $\Omega$, each of size at least $D$.
Then there exists $\Omega'\subseteq\Omega$ such that
\[
|R\cap\Omega'|\geq b
\quad\text{for every }R\in\Rcal,
\qquad
|\Omega'|\leq\frac{32M(b+L)}{D}.
\]
\end{lemma}

\begin{proof}
If $\Rcal=\varnothing$, take $\Omega'=\varnothing$.
We may therefore assume that $\Rcal$ is nonempty,
so $M\geq D>0$.
Select each element of $\Omega$ independently with probability
$p=\min\{1,8(b+L)/D\}$, and let $\Omega'$ be the resulting set.
If $p=1$, then $\Omega'=\Omega$ and $|R\cap\Omega'|\geq D\geq b$
for every $R\in\Rcal$.
Moreover, $D\leq8(b+L)$, so the asserted size bound holds.

Suppose now that $p<1$.
For each $R\in\Rcal$, put $Y_R=|R\cap\Omega'|$ and
$\mu_R=\mathbb E Y_R=p|R|\geq8(b+L)$.
By independence,
\[
\mathbb E(2^{-Y_R})
=(1-p/2)^{|R|}
\leq\exp(-\mu_R/2).
\]
Applying Markov's inequality to $2^{-Y_R}$ gives
\[
\begin{aligned}
\mathbb P(Y_R\leq\mu_R/2)
&\leq2^{\mu_R/2}\mathbb E(2^{-Y_R})\\
&\leq\exp\left(-\frac{1-\log2}{2}\mu_R\right)
\leq\exp(-\mu_R/8),
\end{aligned}
\]
where the last step uses $\log2<3/4$.
Since $b\leq\mu_R/2$, it follows that
$\mathbb P(Y_R<b)\leq\e^{-(b+L)}$.
The union bound therefore yields
\[
\mathbb P\bigl(\exists R\in\Rcal:
|R\cap\Omega'|<b\bigr)
\leq|\Rcal|\e^{-(b+L)}
\leq\e^{-b}
<\frac14.
\]

On the other hand, $\mathbb E|\Omega'|=pM$,
so another application of Markov's inequality gives
$\mathbb P(|\Omega'|>4pM)\leq1/4$.
The probability that either requirement fails is thus
at most $\e^{-b}+1/4<1$.
Consequently, there is a choice of $\Omega'$ satisfying
all the intersection requirements and
$|\Omega'|\leq4pM=32M(b+L)/D$.
\end{proof}

{Fix an intersecting $r$-uniform base family $\Fcal$
with $\tau(\Fcal)=r$, where $r\geq3$, and put
$m=|\Fcal|$ and $\Hcal=\Hcal_{\mathrm{odd}}(\Fcal)$.
We will apply Lemma~\ref{lem:selection} with $\Omega=\Hcal$.}
To control the number of requirements, we count the distinct
sets of missed edges, rather than the vertex sets producing them.

For $T\subseteq X=V(\Fcal)\cup Z$, define
\[
R(T)=\{F\in\Fcal:F\cap T=\varnothing\},
\qquad U(T)=T\cap Z.
\]
The pair $(R(T),U(T))$ determines the members of both
$\Acal$ and $\Bcal$ disjoint from $T$.
The former consist of $R(T)$ together with $Z$ when
$U(T)=\varnothing$; the latter are precisely the sets
$F\cup\{z\}$ with $F\in R(T)$ and $z\in Z\setminus U(T)$.
Since $R(T)\subseteq\Fcal$ and $U(T)\subseteq Z$,
there are at most $2^{m+r}$ possible pairs.

For $S\subseteq V(\Hcal)$, write
${\mathcal M_{\Hcal}(S)}=\{E\in\Hcal:E\cap S=\varnothing\}$.
The three local pairs associated with $S\cap X_i$,
$i=0,1,2$, determine ${\mathcal M_{\Hcal}(S)}$:
an edge of each product family in \eqref{eq:old-4-1}
avoids $S$ precisely when both of its constituent edges do.
Consequently, the family
\[
{\Rcal(\Hcal)}=\{{\mathcal M_{\Hcal}(S)}:
S\in\binom{V(\Hcal)}{2r}\}
\]
satisfies
\begin{equation}\label{eq:old-5-1}
{|\Rcal(\Hcal)|}\leq2^{3(m+r)}\leq\e^{3(m+r)}.
\end{equation}
Moreover, Lemma~\ref{lem:odd-avoid} gives $|R|\geq D$
for every ${R\in\Rcal(\Hcal)}$, where $D=2r(r-1)$.
These estimates allow us to retain sufficiently many edges
in every member of ${\Rcal(\Hcal)}$ while reducing the total number
of edges to quadratic order.

\begin{theorem}\label{thm:odd}
There is an absolute constant $C_2>0$ such that
$m(2r+1)\leq C_2r^2$ for all sufficiently large $r$.
\end{theorem}

\begin{proof}
For sufficiently large $r\geq3$, choose $\Fcal$ as in
Theorem~\ref{thm:kahn}, so that $m=|\Fcal|\leq C_0r$.
{Let $\Hcal=\Hcal_{\mathrm{odd}}(\Fcal)$ be the
cyclic construction in \eqref{eq:old-4-1}, and let
$\Rcal=\Rcal(\Hcal)$ be its family of missed-edge sets
arising from $2r$-sets, as defined above.}
Apply Lemma~\ref{lem:selection} with
\[
M=3m(m+1)r,\qquad
b=2r+1,\qquad
L=3(m+r),\qquad
D=2r(r-1).
\]
The bound on $|\Rcal|$ follows from \eqref{eq:old-5-1},
and every member of $\Rcal$ has size at least $D$.
Also, $D\geq b\geq2$ for $r\geq3$, so all hypotheses
of the selection lemma hold.

We obtain a subfamily $\Hcal'\subseteq\Hcal$ such that
$|R\cap\Hcal'|\geq2r+1$ for every $R\in\Rcal$ and
\begin{equation}\label{eq:old-5-2}
|\Hcal'|
\leq\frac{32M(b+L)}{D}
=\frac{48m(m+1)(3m+5r+1)}{r-1}.
\end{equation}
Regard $\Hcal'$ as a hypergraph on $V(\Hcal)$.
For every $S\in\binom{V(\Hcal)}{2r}$, we then have
\[
\miss_{\Hcal'}(S)
=|{\mathcal M_{\Hcal}(S)}\cap\Hcal'|
\geq2r+1.
\]
Thus $\gamma_{2r}(\Hcal')\geq2r+1$.
As a subfamily of $\Hcal$, the hypergraph $\Hcal'$ remains
intersecting and $(2r+1)$-uniform, and hence is admissible
in the definition of $m(2r+1)$.
Finally, $m\leq C_0r$ makes the right-hand side of
\eqref{eq:old-5-2} at most $C_2r^2$ for an absolute
constant $C_2>0$.
Therefore $m(2r+1)\leq|\Hcal'|\leq C_2r^2$.
\end{proof}

Theorems~\ref{thm:even} and~\ref{thm:odd} give the quadratic
upper bound for all sufficiently large uniformities.
Together with the lower bound in Lemma~\ref{lem:lower},
this proves Theorem~\ref{thm:main}.

Having proved the main theorem, we turn to the extension
to larger avoidance multiplicities.
The odd-uniformity construction provides a quadratic avoidance
guarantee before sparsification. By varying the sampling
probability, we can therefore prescribe a larger avoidance
multiplicity while keeping the number of edges within a constant
factor of the lower bound in Lemma~\ref{lem:lower}.

\begin{corollary}\label{cor:more}
There exist absolute constants $C>0$ and $k_0$ such that,
for every odd integer $k\geq k_0$ and every integer $b$ with
\[
k\leq b\leq\frac{(k-1)(k-3)}2,
\]
we have $kb+1\leq m(k,b)\leq Ckb$.
In particular, $m(k,b)=\Theta(kb)$ uniformly over this range.
\end{corollary}

\begin{proof}
Write $k=2r+1$.
For sufficiently large $r$, choose the base family $\Fcal$
from Theorem~\ref{thm:kahn}, with $m=|\Fcal|\leq C_0r$,
{and form the cyclic construction
$\Hcal=\Hcal_{\mathrm{odd}}(\Fcal)$ in \eqref{eq:old-4-1}.
Let $\Rcal=\Rcal(\Hcal)$ be the corresponding family
of missed-edge sets arising from $2r$-sets.}
It has $M=3m(m+1)r=O(r^3)$ edges.
By Lemma~\ref{lem:odd-avoid}, every $2r$-set misses at least
$D=2r(r-1)$ edges, while \eqref{eq:old-5-1} shows that
these sets determine at most $\e^L$ distinct missed-edge
sets, where $L=3(m+r)=O(r)$.
The assumed range of $b$ is precisely
$2r+1\leq b\leq D$, so Lemma~\ref{lem:selection} applies
{with $\Omega=\Hcal$ and this family $\Rcal$}.

{The selection lemma gives a subfamily
$\Hcal'\subseteq\Hcal$ with $|R\cap\Hcal'|\geq b$
for every $R\in\Rcal$. Regard $\Hcal'$ as a hypergraph
on $V(\Hcal)$. For every $S\in\binom{V(\Hcal)}{2r}$,
we have $\miss_{\Hcal'}(S)
=|\mathcal M_{\Hcal}(S)\cap\Hcal'|\geq b$.
Thus $\gamma_{2r}(\Hcal')\geq b$, and the same selection
lemma bounds its size by}
\[
|\Hcal'|
\leq\frac{32M(b+L)}{D}
=O\bigl(r(b+r)\bigr)
=O(rb)
=O(kb),
\]
where we used $b\geq2r+1$.
All implied constants depend only on the absolute constant
$C_0$ and are therefore independent of $k$ and $b$.
Since $\Hcal'$ remains intersecting and $k$-uniform,
it follows that $m(k,b)\leq Ckb$ for an absolute constant $C$.
The lower bound $m(k,b)\geq kb+1$ follows from
Lemma~\ref{lem:lower}.
\end{proof}

We conclude with the scope of these results and some further questions.
Theorem~\ref{thm:main} determines the order of the Duke--Erd\H{o}s extremal function for every sufficiently large uniformity. The even-uniformity estimate follows from the classical three-block construction and Kahn's theorem; the odd-uniformity construction and its sparsification supply the remaining case. The argument also gives the uniform two-parameter estimate in Corollary~\ref{cor:more}.

The constants in the upper bounds are not optimized. It remains to determine whether $m(k)/k^2$ has a limit, and to improve either the leading lower bound or the constant supplied by the constructions. For larger avoidance multiplicities, it would be useful to obtain $m(k,b)=\Theta(kb)$ over a comparable range for even $k$, and to determine the behaviour beyond the range $b\leq(k-1)(k-3)/2$ established here for odd $k$. These questions concern finer estimates than the order of growth proved in this paper.

\section*{Declaration of generative AI and AI-assisted technologies}
The authors used a large language model during the preparation of this manuscript
for assistance with mathematical exploration, checking arguments,
organization and language editing. The authors are responsible for
the mathematical statements, proofs, references and final text.

\bigskip
{\small
\noindent\textsc{Guoqing Wang}\par
\noindent\textsc{School of Mathematics, Tiangong University}\par
\noindent\emph{Email address:} \texttt{gqwang1979@aliyun.com}\par
\medskip
\noindent\textsc{Linlin Wang}\par
\noindent\textsc{School of Mathematics, China University of Mining and Technology}\par
\noindent\emph{Email address:} \texttt{5400@cumt.edu.cn}\par
\medskip
\noindent\textsc{Jun Yue}\par
\noindent\textsc{School of Mathematics, Tiangong University}\par
\noindent\emph{Email address:} \texttt{yuejun06@126.com}\par
}
\end{document}